\documentclass[12pt]{amsart}
\usepackage{hyperref}
\usepackage{amsfonts}
\usepackage{amsmath,amsthm,amssymb}
\usepackage{amssymb,pdfsync}
\usepackage{graphicx}
\usepackage[mathscr]{eucal}
\usepackage{dsfont}
\usepackage[usenames]{color}
\definecolor{Purple}{rgb}{.7,0.08,0.6}

\usepackage{amssymb,pdfsync}
\usepackage{graphicx}
\usepackage{bm}
\usepackage[mathscr]{eucal}
\usepackage[usenames]{color}
\definecolor{Purple}{rgb}{.7,0.08,0.6}

\usepackage{enumerate}

\theoremstyle{plain}
\newtheorem{Thm}{Theorem}

\newtheorem{Prop}[Thm]{Proposition}

\theoremstyle{definition}

\newcommand{\bwinf}{{\overset{\bm .}{W}}_{\!\infty}}

\newcommand{\Aop}{\mathcal A}

\newcommand{\Cop}{\mathcal C}

\newcommand{\Sop}{\mathcal S}

\renewcommand{\bar}{\overline}
\newcommand{\dee}{\partial}
\newcommand{\deebar}{\overline\partial}
\newcommand{\de}{\mathrm{d}}

\newcommand{\bndry}{b}
\numberwithin{equation}{section}

\begin{document}
\title
[Regularity \& Irregularity]{On Regularity and Irregularity\\
of certain \\
holomorphic singular integral operators}
\author[Lanzani and Stein]{Loredana Lanzani$^*$
and Elias M. Stein$^{**}$}
\thanks{$^*$ Supported in part by the National Science Foundation, award
DMS-1503612.}
\thanks{$^{**}$ Supported in part by the National Science Foundation, award
DMS-1700180.}
\address{
Dept. of Mathematics,       
Syracuse University 
Syracuse, NY 13244-1150 USA}
  \email{llanzani@syr.edu}
\address{
Dept. of Mathematics\\Princeton University 
\\Princeton, NJ   08544-100 USA }
  \thanks{2000 \em{Mathematics Subject Classification:} 30E20, 31A10, 32A26, 32A25, 32A50, 32A55, 42B20, 
46E22, 47B34, 31B10}
\thanks{{\em Keywords}: Hardy space; Cauchy Integral; 
 Cauchy-Szeg\H o projection;  Lebesgue space; pseudoconvex domain; minimal smoothness; measure}
 \thanks{\today}
\begin{abstract} We survey recent work and announce new results concerning 
 two singular integral operators whose kernels are holomorphic functions of the output variable, specifically the Cauchy-Leray integral and the Cauchy-Szeg\H o projection associated to various classes of bounded domains in $\mathbb C^n$ with $n\geq 2$. 
\end{abstract}
\maketitle

\section{Introduction}
\noindent This  is a review of recent and forthcoming work concerning
a menagerie of singular integral operators in several complex variables whose kernels are holomorphic functions of the output variable. (All proofs have appeared or will appear elsewhere.) Our family of operators consists of  the {\em Cauchy-Szeg\H o Projection}, namely 
   the orthogonal projection of  $L^2(\bndry D, \mu)$ onto the holomorphic Hardy space $\mathcal H^2 (\bndry D, \mu)$, as well as various higher-dimension analogs of the Cauchy integral for a planar curve that are collectively known as {\em Cauchy-Fantappi\`e integrals} and include
  the {\em Cauchy-Leray integral} as a particularly relevant example.
 We will henceforth denote such operators $\Sop$ and $\Cop$, respectively.
   Here $D$ is a bounded  domain in complex Euclidean space
$\mathbb C^n$ with $n\geq 2$;
$\bndry D$ is the topological boundary of $D$, while $\mu$ is an appropriate measure supported on $\bndry D$, and we will pay particular attention to two such measures, namely  {\em induced Lebesgue measure} $\Sigma$, 
and the {\em Leray-Levi measure} $\lambda$. 
\vskip0.1in
To be precise, we are interested in {\em the $L^p$-regularity problem for $\Cop$ and for $\Sop$}, that is:
\vskip0.1in
 Determine regularity and geometric conditions on the ambient domain $D$ that grant that
\vskip0.1in
\begin{itemize}
\item[$\bullet$] $\Cop : L^p(\bndry D, \lambda)\to L^p(\bndry D, \lambda)$
is bounded for all $p\in (1, \infty)$.
\vskip0.1in
\item[$\bullet$] $\Sop : L^p(\bndry D, \Sigma)\to L^p(\bndry D, \Sigma)$
is bounded for $p$ in an interval of maximal size about $p_0=2$.
\end{itemize}
In complex dimension 1 (that is, for $D\Subset \mathbb C$) these problems are well understood; here we focus on dimension $n\geq 2$. 

\vskip0.1in

We point out that the Cauchy-Leray integral is a Calder\'on-Zygmund operator, thus $L^p$-regularity for $p=2$ is equivalent to 
regularity in $L^p$ for $1<p<\infty$. On the other hand,
the Cauchy-Szeg\H o projection 
  is automatically bounded in $L^2(\bndry D, \mu)$ but $L^2$-regularity does not guarantee $L^p$-regularity  for $p\neq 2$: indeed, establishing $L^p$-regularity of $\Sop$ for $p\neq 2$ is, in general, a very difficult problem.
The main difficulty stems from the fact that the Schwartz kernel for $\Sop$ (that is the {\em Cauchy-Szeg\H o kernel}\,) is almost never explicitly available, even in the favorable setting when $D$ is smooth and strongly pseudoconvex,  so direct estimates cannot be performed and one has to rely on other methods, such as
 asymptotic formulas analogous to those obtained by C. Fefferman \cite{F} (for the Bergman kernel)  and Boutet de Monvel-Sj\"ostrand \cite{BS},
or a paradigm discovered
 by N. Kerzman and E. M. Stein \cite{KeSt-1} that relates $\Sop$ to a certain Cauchy-Fantappi\`e
 integral associated to $D$ (the {\em Kerzman-Stein identity}).

About 30 years later, a surge of interest 
in 
singular integral operators in a variety of ``non-smooth'' settings led us to a new examination of these problems from the following 
  point of view:
\vskip0.03in

 {\em to what extent is the $L^p$-boundedness of the aforementioned operators reliant upon the boundary regularity and (natural to this context) 
upon the amount of convexity of the ambient domain $D$?} 
\vskip0.03in

\noindent  Stripping away the smoothness assumptions  brings  to the fore 
     the geometric interplay between the operators and the domains on which they act:
 it soon became apparent that new ideas and techniques were needed, even to deal with rather tame singularities such as the class $C^{2, \alpha}$. 
 The following results were proved in \cite{LaSt-8} and \cite{LaSt-5}.
 \vskip0.1in
 \begin{itemize}
\item[(I)] $\Sop$: $L^p(\bndry D, \Sigma)\to L^p(\bndry D, \Sigma)$
is bounded for {$1<p<\infty$} if
\vskip0.1in
\begin{itemize}
\item[{\tt (i.)}] $D\Subset \mathbb C^n$ is strongly pseudoconvex, and
\item[{\tt (ii.)}] $\bndry D$ is of class $C^2$.
\end{itemize}
\vskip0.2in
\item[(II)] $\Cop$: $L^p(\bndry D, \lambda)\to L^p(\bndry D, \lambda)$
is bounded for {$1<p<\infty$} if
\vskip0.1in
\begin{itemize}
\item[{\tt (i.)}] $D\Subset \mathbb C^n$ is strongly $\mathbb C$-linearly convex, and
\item[{\tt (ii.)}] $\bndry D$ is of class $C^{1,1}$.
\end{itemize}
\end{itemize}
\vskip0.004in
 The purpose of this note is to summarize the main points in the proof of (I) and (II), and to announce   new results 
 that will appear in forthcoming papers \cite{LaSt-6} and \cite{LaSt-10}
 pertaining
  the optimality of the assumptions made in (I) and (II).

\section{The $L^p$-regularity of the Cauchy-Leray integral}

 Here we
  we take as our model the seminal one-dimensional theory 
of Calder\`on \cite{Ca},
 Coifman-McIntosh-Meyer 
\cite{CMM}, and David 
\cite{Da} for the Cauchy integral of a planar curve, and in particular its key theorem on the {\em Lipschitz} case. 
As is well known, the initial result was the classical theorem of M. Riesz for the Cauchy integral on the unit disc (i.e. the Hilbert transform on the circle); 
the standard proofs which developed from this then allowed an extension to a corresponding result where the disc is replaced by a domain $D\subset\mathbb C$ whose boundary is relatively smooth, i.e. of class $C^{1, \alpha}$, for $\alpha >0$. However, going beyond that to the limiting case of regularity, namely 
$C^1$ and other variants ``near $C^1$'', required further ideas. The techniques introduced in this connection led to significant developments in harmonic analysis such as the ``$T(1)$ theorem'' and various aspects of multilinear analysis and analytic capacity, \cite{Ch-1}, \cite{MeCo}, \cite{To-1}, \cite{To-2}. The importance of those advances suggests the following fundamental question: {\em what might be the corresponding results for the Cauchy integral in several variables}.  However, in the context of higher dimension 
 geometric obstructions arise (pseudoconvexity or, equivalently, lack of conformal mapping) which in the one-dimensional 
  setting
  are irrelevant. As a consequence,  there is no canonical notion of {\em holomorphic} Cauchy kernel: all such kernels must be domain-specific. Indeed, the only  kernel that can be deemed ``canonical''\footnote{``canonical'' in the sense that it is
  the restriction to $\bndry D$ of a {\em universal} kernel defined in $\mathbb C^n\times \mathbb C^n\setminus\{z=w\}$.}
   is the Bochner-Martinelli kernel \cite{Ky},
   but such kernel is nowhere holomorphic and thus of no use in the applications described below. 
  One is therefore charged with the further
  task of   {\em constructing} a holomorphic kernel that is fitted to the specific geometry of the domain  and, after that, with supplying proof of
    regularity of the resulting singular integral operator.
 As in the one-dimensional setting,  this theory was first conceived within the context of  {\em smooth} ambient domains; if the domain is
  not sufficiently smooth (of class $C^{2,\alpha}$ or better)
    the original kernel constructions by Henkin and Ramirez,
     \cite{He} and
      \cite{Ram}, and  the ``osculation by the Heisenberg group'' technique
       in Kerzman-Stein
       \cite{KeSt-1} are no longer applicable. 
     In \cite{LaSt-5} it is shown that
the  $T(1)$-theorem technique for a space of homogeneous type fitted to the geometry and regularity of the ambient domain can be applied to prove $L^p(\bndry D, \mu)$-regularity, for $1<p<\infty$ of the {\em Cauchy-Leray integral}\,:
   \begin{equation}\label{E:Leray-Int}
  \Cop f(z) =\frac{1}{(2\pi i)^n}\int\limits_{w\in\bndry D}
  f(w)\frac{\Omega (z, w)\wedge (\de_w \Omega (z, w))^{n-1}}
  {\langle \Omega (z, w), w-z\rangle^n},\ \ \ z\in D
  \end{equation}
whenever $D\subset\mathbb C^n$ is  a bounded, strongly $\mathbb C$-linearly convex domain  
whose boundary satisfies the minimal regularity condition given by the class $C^{1, 1}$ (that is, the domain admits a defining function $\rho$ of class $C^{1,1}$).  Here the {\it generating 1-form} $\Omega (z, w)$ is 
the complex gradient of the domain's defining function (and we should point out that the definition of $\Cop$ is independent of the choice of defining function $\rho$).
More precisely: $\Omega (z, w)$ = $j^*\dee\rho (w)$, where $j^*$ denotes the pull-back under the inclusion map $j:\, \bndry D\,\hookrightarrow\, \mathbb C^n$. 
    The boundary measure $\mu$ belongs to a family that includes induced Lebesgue measure $\Sigma$, as well as the Leray-Levi measure
\begin{equation}\label{E:Leray-meas}
       \de\lambda(w) := (2\pi i)^{-2n}
       j^*\left(\dee\rho\wedge(\deebar\dee\rho)^{n-1}\,\right)(w)\, ,\quad w\in\bndry D.
       \end{equation}
       \vskip-0.05in
We remark that under our assumptions (class $C^{1,1}$) the factor
   $\deebar\dee\rho$ in the definition of the Leray-Levi measure $\lambda$, as well as  the factor $\de_w \Omega (z, w)$ in the
   Schwartz kernel for $\Cop$, are only in $L^\infty(\mathbb C^n)$ and therefore may be undefined on $\bndry D$ because  the latter is a zero-measure
   subset  of $\mathbb C^n$, however it turns out that the {\em tangential component} of each of 
   $\deebar\dee\rho(w)$ and $\de_w \Omega (z, w)$, namely $j^*\deebar\dee\rho(w)$ and $j^*\de_w \Omega (z, w)$, are in fact meaningful, leading to a  kernel that is  well-defined even in our singular context.

\section{Counter-examples
to the $L^p$-theory for the Cauchy-Leray integral}
 In  \cite{LaSt-6} we construct two examples that establish the optimality of the assumptions made on the ambient domain for the Cauchy-Leray integral $\Cop$. Both examples are real ellipsoids of the form
\begin{equation}\label{E:exCL}
D_{r, q} \ :=\ \bigg{\{}(z_1, z_2)\in\mathbb C^2\ \ \bigg|\ \   \left|\mathrm{Re}\, z_1\right|^r\, +\, 
\left|\mathrm{Im}\, z_1\right|^q\, +\, |z_2|^2\, < \, 1\bigg{\}}\, .
\end{equation}
For the first example, $(r, q)= (2, 4)$; in this case the domain is smooth, strongly pseudoconvex and strictly convex, but it is not strongly $\mathbb C$-linearly convex. In the second example, $(r, q)= (m, 2)$ for any $1<m<2$; this domain is strongly $\mathbb C$-linearly convex but is only of class $C^{1, m-1}$ (and no better). In both cases we show that the associated Cauchy-Leray integral $\Cop$ is well-defined on a dense subset of $L^p(\bndry D_{r, q}, \mu)$ but {\em does not extend to a bounded operator}: 
$L^p(\bndry D_{r, q}, \mu) \mapsto L^p(\bndry D_{r, q}, \mu)$ for any $1<p<\infty$.
Specifically, we prove that there is a function $f\in C^1(\bndry D_{r, q})$  supported in a proper subset of $\bndry D_{r, q}$  
 such that  
  \begin{itemize}
  \vskip-0.3in
  \item[{\tt(i.)}] $\Cop f (z)$ can be defined as an absolutely convergent integral whenever $z\in \bndry D_{r, q}$  is at positive distance from the support of $f$. 
  \vskip0.05in
\item[{\tt (ii.)}]
The inequality:
$ \|\Cop (f)\|_{L^p(S,\, d\mu)}\leq A_p\|f\|_{L^p (\bndry D_{r, q},\, d\mu)}$
(with $A_p$ independent of $f$ and $S$) {\em fails} 
whenever
 $S\subset\bndry D_{r, q}$
 is disjoint from the support of $f$. 
 \end{itemize}
Here $\mu$ is a boundary measure that belongs to a family that includes induced-Lebesgue measure $\Sigma$ and the Leray-Levi measure $\lambda$, as well as Fefferman's  measure $\mu^F$, see \cite{F}; for the first example (the smooth domain $D_{2, 4}$) all such measures are mutually absolutely continuous.  For the second example (the non-smooth domain $D_{m, 2}$) these measures are essentially different, yet the counter-example holds in all cases.

The main tool for proving {\tt (i.)} and {\tt (ii.)} is a scaling and limiting process that transfers the problem to  specific,  unbounded smooth domains, namely $\{ \, 2\,\mathrm{Im}\, z_2> (\mathrm{Re}\, z_1)^2\}$ in the first case, 
and $\{ \, 2\,\mathrm{Im}\, z_2> \left|\mathrm{Re}\, z_1\right|^m\}$ in the second. On the unbounded domains, explicit computations are carried out to prove failure of the $L^p$-boundedness of the transported operator. 
\noindent There is also the matter of showing that the Cauchy-Leray integral
for $D_{r, q}$ maps $L^p$ into the holomorphic Hardy space $\mathcal H^p$: this question
is addressed in \cite{LaSt-7} where it is shown that
\begin{itemize}
\item[{\tt (iii.)}] $\Cop f (z)$, for $z\in\bndry D_{r, q}$ as in item {\tt(i.)} above, arises as ``boundary value'' of a function $F$ holomorphic in $D_{r, q}$.
\end{itemize}
The proof of {\tt (iii.)} requires three different approaches, each tailored to the particular type of singularity displayed by the example under consideration: in dealing with 
the non-smooth domain $D_{m, 2}$ one has to distinguish the case when $1<m\leq 3/2$ from the case $3/2<m<2$: in the second case, a global integration by parts gives that 
$\Cop f (z)$ is the restriction to $\bndry D_{m, 2}$ of a holomorphic $F\in C^1(\bar D_{m, 2})$ . On the other hand, when $1<m\leq 3/2$ such method is no longer viable but  we show nonetheless, that $\Cop f$ extends to a holomorphic   $F$ that is continuous everywhere on $\bar D_{m, 2}$ except for a 0-measure subset of the boundary (namely the sphere $\{ |\text{Re}\, z_1|^2 + |z_2|^2 =1\}$).

 Finally,  the lack of strong $\mathbb C$-linear convexity in the first example (the domain $D_{2, 4}$)  prevents us from  carrying a global integration by parts: instead, one shows that $\Cop f$ extends to a holomorphic $F\in C(\bar{D}_{2, 4})$ by using a local integration by parts which depends on the location of the coordinate patch with respect to the ``flat'' part of the boundary. 
It should be noted that an earlier result in 
Barrett-Lanzani 
\cite{BaLa}  already gave an example with irregularity in $L^2(\bndry D, \mu)$, however 
the less explicit and more complex nature of the construction did not provide insight for $L^p(\bndry D, \mu)$ when $p\neq 2$.

\section{The $L^p$-regularity of the Cauchy-Szeg\H o projection}

\subsection{Discussion of the problem}

        We recall that the Cauchy-Szeg\H o projection $\Sop$ is the unique, {\em orthogonal} (equivalently, {\em selfadjoint}) projection operator of 
 $L^2(bD, \Sigma)$ onto the Hardy space of holomorphic functions; here $\Sigma$ is the induced  Lebesgue measure on $\bndry D$. As mentioned earlier, one must come to terms with the fact that, in general,  orthogonal projections are {\em not} Calder\`on-Zygmund operators,
 thus  $L^p$-regularity for $p\neq 2$
 does not follow from $L^2$-regularity; also, one may have $L^p$-regularity only for $p$ in a proper {\em sub-interval} of $(1, \infty)$, see e.g. \cite{LaSt-1}.
  (By contrast, for a Calder\`on-Zygmund operator boundedness in $L^2$ {\em implies} boundedness in $L^p$
  for $1<p<\infty$.)
  Regularity properties of the Cauchy-Szeg\H o  projection, in particular $L^p$-regularity, have been the object of considerable interest for more that 40 years.
        When the boundary of the domain $D$ is sufficiently smooth, decisive results were obtained in the following settings:
{\em (a)}, when $D$ is strongly pseudoconvex \cite{BS}, \cite{PS};
{\em (b)}, when $D \subset \mathbb C^2$ and its boundary is of finite type \cite{Mc-1}, \cite{NRSW};
{\em (c)}, when $D \subset \mathbb C^n$ is convex and its boundary is of finite type \cite{Mc-1}, \cite{McSt}; 
 {\em (d)}, when 
 $D \subset \mathbb C^n$ is of finite type and its Levi form is diagonalizable \cite{CD}. Related results include \cite{AS-1}-\cite{AS-3}, \cite{Ba}, \cite{BoLo}, \cite{KrPe}, \cite{NaPr}, \cite{PoSt}, \cite{Ro-1}, \cite{Ro-2}, \cite{Z}.
 The main difference when dealing with the situation when $D$ has lower (in fact minimal) regularity than 
 the setting of the more regular domains treated in {\em (a) -- (d)}, is 
 that in each of those cases known formulas for the Cauchy-Szeg\H o kernel, or at least size estimates, played a decisive role. In our general situation such estimates are unavailable and one must proceed by a different analysis that relies upon {\tt (i.)}, the $T(1)$-theorem technique of \cite{LaSt-5} and {\tt (ii.)}, a new, tricky variant of the original Kerzman-Stein paradigm \cite{KeSt-1} described below.

\subsection{$L^p$-regularity of the Cauchy-Szeg\H o projection.}  Strong $\mathbb C$-linear convexity
              {\it implies} strong pseudoconvexity whenever
        the domain enjoys enough regularity for the latter to be meaningful.   
        In  \cite{LaSt-8} 
        some of
       the techniques from
        \cite{LaSt-5} are adapted to study the $L^p$-regularity problem for the Cauchy-Szeg\H o  projection
       of strongly Levi-pseudoconvex domains $D\Subset \mathbb C^n$ with minimal boundary regularity, namely the class $C^2$  (which is the minimal regularity for strong Levi-pseudoconvexity to hold), 
 leading to the conclusion  that
$L^p$-boundedness of $\Sop$ holds in the full range $1<p<\infty$.
 As mentioned above, in this general setting a direct analysis of the
 Cauchy-Szeg\H o kernel does not lead to the desired result. Instead, our starting point is the original Kerzman-Stein paradigm  
  \cite{KeSt-1}
   for domains that are sufficiently smooth: this proceeded by constructing
    a holomorphic Cauchy-Fantappi\`e integral $\Cop$ in the same spirit of \eqref{E:Leray-Int} but for a different choice of generating form $\Omega$. 
The 
 analysis of
  $\Sop$ 
   begins with the representation:
$\Cop=\Sop(I-\mathcal A)$ on $L^2(\bndry D, \Sigma)$, where 
$I$ is the identity and
$\mathcal A$ denotes the difference of $\Cop$ and its formal
$L^2$-adjoint, that is: $\mathcal A= \Cop^*-\Cop$.
This identity follows from the
fact that, just like the Cauchy-Szeg\H o projection, the Cauchy-Fantappi\`e integral $\Cop$ is also a projection  of $L^2(\bndry D, \Sigma)$ onto the
  holomorphic Hardy space\footnote{It is failure of this property that renders the Bochner-Martinelli integral unsuitable for the analysis of $\Sop$.} (albeit not the orthogonal projection!). \nopagebreak  In particular, since $\mathcal A^* =-\mathcal A$ it follows that the operator $(I-\mathcal A)$ is invertible in $L^2(\bndry D, \Sigma)$ with bounded inverse, and we obtain:
\begin{equation}\label{E:KS}
\Sop=\Cop(I-\mathcal A)^{-1}\quad \text{on}\  L^2(\bndry D, \Sigma).
\end{equation}

  Kerzman and
  Stein
\cite{KeSt-1} proved
that if the (strongly pseudoconvex) domain is sufficiently smooth (e.g. of class
$C^3$)
  the singularities of $\Cop$ and $\Cop^*$ cancel out and as a result
 $\mathcal A$ is ``small'' in the sense that it
is compact in $L^2(\bndry D, \Sigma)$ (indeed smoothing); from this it follows that the righthand side of 
 the above identity
 is bounded in
$L^p(\bndry D, \Sigma)$ for all $1<p<\infty$ and therefore so is $\Sop$, giving the solution to the $L^p$-regularity problem for $\Sop$ in the full range $1<p<\infty$. 

If the domain is only of class $C^2$ this argument
is no longer applicable because $\Aop$ in general fails to be compact on $L^2(\bndry D, \Sigma)$, see \cite{BaLa}. Instead, in
\cite{LaSt-8} 
we work
 with a {\em family} of holomorphic Cauchy-Fantappi\`e integrals
   $\{\Cop_\epsilon\}_\epsilon$ whose kernels are constructed via a first-order perturbation of the
   Cauchy-Leray kernel \eqref{E:Leray-Int}
  that makes use of a smooth approximation $\{\tau_\epsilon\}_\epsilon$ of certain second-order derivatives of the
   defining function of the domain. 
   As in the case of the Cauchy-Leray integral $\Cop$, here there
   are two boundary measures at play: the induced Lebesgue measure $\Sigma$, and the Leray-Levi measure  $\lambda$, see \eqref{E:Leray-meas},
 which   in this new context is absolutely continuous with respect to $\Sigma$ because of the relation
 \begin{equation}\label{E:Leray-surf}
     \de\lambda(w)\ \approx\ |\varphi (w)|\de\Sigma (w)\, , \ \ w\in\bndry D\, 
      \end{equation}
  where $\varphi (w)$ is the determinant of the Levi matrix. The operators $\{\Cop_\epsilon\}_\epsilon$ are then seen to be bounded in 
  $L^p(\bndry D, \lambda)$ {\em and} $L^p(\bndry D, \Sigma)$ for all $1<p<\infty$ by an application of the  $T(1)$-theorem.
  On the other hand, in defining the Cauchy-Szeg\H o projection it is imperative to specify the underlying
measure for $\bndry D$ that arises in the notion of orthogonality
that is being used. Correspondingly, we now have two distinct Cauchy-Szeg\H o projections $\Sop_\Sigma$ and $\Sop_\lambda$ but these, in our general setting, are not directly related to one another. 
 It turns out that
 the Leray-Levi measure $\lambda$ has a ``mitigating'' effect 
  that leads to a new  smallness argument for the difference
   $\Cop_\epsilon^\dagger - \Cop_\epsilon$ that occurs {\em when the adjoint $\Cop_\epsilon^\dagger$ is computed with respect to $\lambda$}.
 While the $\{\Cop_\epsilon\}_\epsilon$  do not approximate $\Sop_\lambda$
   (in fact the norms of the $\Cop_\epsilon$ are in general unbounded as $\epsilon\to 0$),
 we show that for each fixed $1<p<\infty$ (in fact for $p<2$) there is $\epsilon = \epsilon (p)$
 such that $\Cop_\epsilon^\dagger-\Cop_\epsilon$ splits as the sum $\mathcal B_\epsilon +\Aop_\epsilon$, where 
 $\mathcal B_\epsilon:
  L^p (\bndry D, \lambda)\to C(\bndry D)$, and $\|\Aop_\epsilon\|_{L^p\to L^p}\leq \epsilon$: this is the new, ``tricky'' variant of the original Kerzman-Stein paradigm
   that was alluded to earlier, and it
     gives us
the identity  
 \begin{equation}\label{E:KS-new}
 \Sop_\lambda = (\Sop_\lambda\mathcal B_\epsilon+\Cop_\epsilon)(I-\mathcal A_\epsilon)^{-1}\quad 
 \text{ in}\quad L^2(\bndry D, \lambda).
 \end{equation}

Then one proves that the righthand side is bounded on $L^p(\bndry D, \lambda)$ (here we also use that $p<2$ and that $D$ is bounded) and we conclude that  $\Sop_\lambda$ is bounded
in $L^p(\bndry D, \lambda)$ whenever $1<p<2$; the result for $p>2$  follows by duality. 
    A similar argument is needed to
   treat $\Sop_\Sigma$, but there is no direct way to show smallness for $\Cop_\epsilon^* - \Cop_\epsilon$ when the adjoint $\Cop_\epsilon^*$ is computed with respect to the induced Lebesgue surface measure 
  $\Sigma$. Instead, one recovers such smallness from the corresponding result for 
 $\Cop_\epsilon^\dagger - \Cop_\epsilon$, by observing that
 $\Cop_\epsilon - \Cop_\epsilon^* =
  \Cop_\epsilon - \Cop_\epsilon^\dagger +|\varphi|^{-1}\big[|\varphi|, \Cop_\epsilon^\dagger\big]$, 
  where $\varphi$ is as in \eqref{E:Leray-surf}, and by 
  controlling the size of the operator norm of
   the commutator $\big[|\varphi|, \Cop_\epsilon^\dagger\big]$.
 \vskip0.05in

To complete the proof one also needs the requisite
   representation formulae and density results for the holomorphic 
   Hardy spaces 
    of the domains that satisfy the minimal boundary regularity conditions
    stated in (I) and (II): these are obtained in \cite{LaSt-9}.
    
    \section{A counter-example
to the $L^p$-theory for the Cauchy-Szeg\H o projection}
 The forthcoming work \cite{LaSt-10} investigates a long-standing open question concerning 
  {\em $L^p$-irregularity}
   of the Cauchy-Szeg\H o projection for the {\em Diederich-Forn\ae ss worm domains}:
   \vskip-0.2in
 \begin{equation}\label{E:worms}
W_{k, h}:=\left\{(z_1, z_2)\in\mathbb C^2\, ,\ \left|z_2-ie^{ih(|z_1|)}\right|^2<1-k(|z_1|)\right\}\, .
\end{equation}
\vskip-0.05in
\noindent  Appropriate choices of the functions $h$ and $k$ produce  
  domains  that are
  smooth and pseudoconvex but only {\em weakly} pseudoconvex along a 2-dim subset of their topological boundary.  (The nick-name ``worm'' is meant to illustrate winding caused by the argument $h(|z_1|)$.)

Developed by Diederich and Forn\ae ss in 1977 as  examples 
 of smooth, weakly pseudoconvx domains with non-trivial Nebenh\"ulle\footnote{ the domain is pseudoconvex but cannot be ``exhausted'' by smooth pseudoconvex ``super-domains''.}, the class \eqref{E:worms} has
  since proved to be a reliable source 
of counter-examples to a variety of phenomena in complex function theory. 
Of special relevance here 
are the seminal paper \cite{Ba}
and the related work \cite{KrPe} that prove {\em irregularity} of the Bergman projection\footnote{that is, the orthogonal projection of $L^2(D, dV)$ onto the Bergman space $A^2(D) := \vartheta(D)\cap L^2(D, dV)$.} for the worm domain in the Sobolev- and Lebesgue-space scales, respectively, when the following choices are made for $h$ and $k$:
\vskip-0.2in
\begin{equation}\label{E:DFworm}
h(|z_1|):= \log\!|z_1|^2\, ;\qquad k(z_1):= \phi (h(|z_1|))
\end{equation} 
\vskip-0.04in
\noindent with $\phi$  a smooth, non-negative even function chosen so that $W_{h, k}$ is smooth, bounded, connected and pseudoconvex, and moreover $\phi^{-1}(0) = \{\, |t|\, \leq \beta-\pi/2\, \}$ for fixed, given $\beta>\pi/2$.

In contrast with the situation for the Cauchy-Leray integral, 
the Cauchy-Szeg\H o and Bergman projections are {\em always} bounded in $L^2$ (that is for $p=2$)
so in this context ``$L^p$-irregularity'' should be interpreted as
 ``failure of $L^p$-regularity in the {\em full range $1<p<\infty$}''. 
 
 The results described henceforth will appear in \cite{LaSt-10}.
  \begin{Thm}[Main result]\label{T:main}
  For any $p\neq 2$ there is $\beta = \beta(p)>\pi/2$ such that for 
  $W =W_{h, k}$ with $h$, $k$ as in \eqref{E:DFworm}, the Cauchy-Szeg\H o projection
 associated to $W$ is not bounded: $L^p(\bndry W, \Sigma)\to L^p(\bndry W, \Sigma)$.
 \end{Thm}
 Here $\Sigma$ is induced Lebesgue measure for $\bndry W$. The strategy of proof is similar in spirit to 
the original arguments \cite{Ba} and \cite{KrPe} for the Bergman projection (which also 
 inspired the strategy of proof for the examples for the Cauchy-Leray integral described in the previous section):
one starts with a (biholomorphic) scaling of the original domain $W$ leading to a family of smooth domains $\{W_\lambda\}_\lambda $;  then a limiting process
 transfers the $L^p$-regularity problem
 to a specific, unbounded limiting domain $W_{\!\infty}$.  On the latter, explicit computations are carried out that prove failure of $L^p$-regularity of the relevant operator for $W_{\!\infty}$. The scaling and limiting arguments then allow to percolate failure of $L^p$-regularity back to $W$ via a suitable transformation law under the scaling map. 
 
 When carrying out this scheme for the Cauchy-Szeg\H o projection  several new obstacles arise that were non existent
 in the analysis of the Bergman projection and of the Cauchy-Leray integral: here we focus on just one, namely the fact that
   the limiting domain $W_{\!\infty}$ is unbounded and {\em non-smooth} (it is a Lipschitz domain),  thus  for $W_{\!\infty}$
 there is no canonical notion of holomorphic Hardy space nor of Cauchy-Szeg\H o projection
  (by contrast, the definition of the Bergman space $A^2(W_{\!\infty}, dV)$ is standard, and so is the associated Bergman projection).
 It is not hard to see that the topological boundary of
   $W_{\!\infty}$
  splits into three distinct parts: two of these, denoted $\bwinf$ and ${\overset{\bm.\bm.}{W}}_{\!\infty}$ have full induced-Lebesgue  measure, while the third part is the distinguished boundary $d_bW_{\!\infty}$. In 
  \cite{MoPe-1} the authors
 prove irregularity of the
 Cauchy-Szeg\H o projection associated to $d_bW_{\!\infty}$ (defined with respect to induced Lebesgue measure
  for $d_bW_{\!\infty}$). However the small size of  the distinguished boundary (it is a codimension-1 subset of the topological boundary) makes it impossible to percolate the result for $d_bW_{\!\infty}$
   back to the Cauchy-Szeg\H o projection for  the full boundary of the original worm $W$.
    Here
   we focus instead on the full-measured part of the boundary 
    denoted 
  $\bwinf$ because this particular piece of the boundary
  supports a natural notion of ``quasi-product measure'' $\mu_\infty$ that 
   captures the main features of the full boundary of $W_{\!\infty}$, as indicated by the following key observation:
  \begin{Prop}
 Suppose that $\displaystyle{F\in C_0\! \left(\overline{\underset{\lambda >0}{\cup}W_\lambda}\right)}$. Then
$$\displaystyle{
\lim\limits_{\lambda \to \infty} \int\limits_{\bndry W_\lambda}\!\!\! F\, d\mu_\lambda =
\int\limits_{\bwinf}\!\!\! F\, d\mu_\infty\, .
}$$
\end{Prop}
Here $\mu_\lambda$ is the transported induced Lebesgue measure for $W$ via the scaling map.  
(In fact a more sophisticated version of the above proposition is needed, one that is valid for $F$ in a larger function space that is dense
in 
the Hardy space for $W$,
 but the above already provides the required ``supporting evidence''.)

 It turns out that the quasi-product measure $\mu_\infty$ leads to a meaningful notion of Hardy space for $\bwinf$ and furthermore, that
the topological boundary of the original (smooth) worm $W$ also supports a ``quasi-product'' measure $\mu_0$ that is mutually absolutely continuous with respect to induced Lebesgue measure $\Sigma$ 
 and enjoys a certain stability under the scaling maps, leading us to the following result:

 \begin{Thm}
 Let $\Sop_\infty$ denote the Cauchy-Szeg\H o projection for $H^2(\bwinf, \mu_\infty)$, and let 
 $\Sop_{\bndry W}$ denote the Cauchy-Szeg\H o projection for $H^2(\bndry W, \Sigma)$.
 \vskip0.1in
If\ \ $\Sop_{\bndry W}: L^p(\bndry W, \Sigma) \to L^p(\bndry W, \Sigma)$
is bounded, then

$\Sop_\infty: L^p(\bwinf, \mu_\infty) \to L^p(\bwinf, \mu_\infty)\quad \text{is bounded and}\ \
$
$$
\| \Sop_\infty\|_{L^p(\bwinf, \,\mu_\infty)\circlearrowright}\ \leq\ 
\| \Sop_{\bndry W} \|_{L^p(\bndry W, \,\Sigma)\circlearrowright}\, .
$$
\end{Thm}

\noindent Finally, a direct examination shows that $\Sop_\infty$ is {\em unbounded} on $L^p(\bwinf, \mu_\infty)$, giving us the proof of Theorem \ref{T:main}.

\vskip0.02in

\vskip0.2in
\begin{center}
{\sf Bibliography}
\end{center}

\end{document}